\documentclass{amsart}

\usepackage{amsmath}
\usepackage{amssymb}
\usepackage{amsfonts}

\usepackage{amsmath}
\usepackage{amsfonts}
\usepackage{amssymb}
\usepackage{amsthm}
\usepackage{hyperref}
\usepackage{graphicx}
\usepackage{subcaption} 
\usepackage{xcolor}

\usepackage{euscript}
\usepackage{setspace}
\usepackage{tikz}
\usetikzlibrary{calc, intersections, through, backgrounds, positioning}

\newcommand{\Z}{\mathbb{Z}}

\newcommand{\C}{\mathbb{C}}
\newtheorem{theorem}{Theorem}

\newtheorem{lemma}{Lemma}

\theoremstyle{remark}

\newcommand{\R}{\mathbb R}

\begin{document}
\title{Dynamical versions of Hardy's  uncertainty principle: a survey
}

\author{
	Aingeru Fern\'andez-Bertolin}
	\address{AFB: Universidad del Pa\'is Vasco /Euskal Herriko Unibertsitatea (UPV/EHU), Dpto. de Matem\'aticas, Apartado 644, 48080 Bilbao, Spain}
	\email{aingeru.fernandez@ehu.eus}
	\author{Eugenia Malinnikova}	
	\address{EM: Department of Mathematics, Stanford University, Stanford, CA, USA \newline
	Department of Mathematical Sciences, Norwegian University of Science and Technology, Trondheim, Norway}
	\email{eugeniam@stanford.edu}
	
	\thanks{EM was partially supported by NSF grant DMS-1956294 and by the Research Council of Norway, project 275113. AFB was partially supported by ERCEA Advanced Grant 2014 669689 - HADE, by the project PGC2018-094528-B-I00 (AEI/FEDER, UE) and acronym ``IHAIP'', and by the Basque Government through the project IT1247-19.}
	
\begin{abstract} The Hardy uncertainty principle says that no function is better localized together with its Fourier transform than the Gaussian. The textbook proof of the result,  as well as one of the original proofs by Hardy, refers to the Phragm\'en-Lindel\"of theorem. In this note we first describe the connection of the Hardy uncertainty to the Schr\"odinger equation, and give a new proof of Hardy's result which is based on this connection and the Liouville theorem. The proof is related to the second proof of Hardy, which has been undeservedly forgotten. Then we survey the recent results on  dynamical versions of Hardy's theorem.
\end{abstract}

	\maketitle

\section{Introduction}
There are many mathematical interpretations of the uncertainty principle, which states that the position and momentum of a quantum particle cannot be measured simultaneously, or that a signal cannot be well-localized  both in time and in frequency. All of them refer to a double representation of a function, classically this is the function itself and its Fourier transform, though more recent versions of the uncertainty principle use some form of joint time-frequency representation, for example the short-time Fourier transform. Each uncertainty principle has an interesting and developing story, in this note we tell only one of them.

The most famous uncertainty principle was introduced by  Werner Heisenberg in 1927, and  its mathematical  formulation was given by  Earle Hesse Kennard and Hermann Weyl shortly after. It says that 
\begin{equation}\label{eq:Hei}
\int_{\R^d}|x|^2|f(x)|^2\int_{\R^d}|\xi|^2|\widehat{f}(\xi)|^2\ge
 \frac{d^2}{4}\|f\|_2^4
\end{equation}
 for all $f\in L^2(\R^d),$ or
equivalently,
\[\int_{\R^d}|x|^2|f(x)|^2+\int_{\R^d}|\xi|^2|\widehat{f}(\xi)|^2\ge d\int_{\R^d}|f|^2.\]
We always use the following normalization of the Fourier transform on $\R^d$,
\[\widehat{f}(\xi)=\frac{1}{(2\pi)^{d/2}}\int_{\R^d}f(x)e^{-ix\cdot\xi}dx.\]
It is well-known that the Fourier transform is an isometry of  $L^2(\R^d)$.

The equality in Heisenberg's uncertainty principle \eqref{eq:Hei} is attained when $f$ is a generalized Gaussian function, i.e., $f(x)=\exp(-(Ax,x))$, where $A$ is a positive definite matrix. The fact that the Gaussian is the best localized function in time and frequency was also recognized by English mathematician Godfrey H. Hardy in 1933, in the formulation of the uncertainty principle that now bears his name. Hardy attributed the remark that a function and its Fourier transform "cannot be very small" to Norbert Wiener and proved the following one dimensional result 

\begin{theorem}\label{th:H}
Let $f\in L^2(\R)$ satisfy $|f(x)|\le Ce^{-a|x|^2}$ and $|\widehat{f}(\xi)|\le Ce^{-b|\xi|^2}.$ If $ab>1/4$ then $f=0$ and if $ab=1/4$ then $f(x)=ce^{-a|x|^2}.$
\end{theorem}

  In his original article \cite{H}, Hardy gave  two different proofs, both refer to holomorphic functions and use some results of complex analysis. The first one employs the Phragm\'en-Lindel\"of principle for entire functions. This proof or its variations can be found in many textbooks, see for example \cite{HJ, Tan, S}. The second one also refers to entire functions, but makes use of the Liouville theorem  only (at least for the case when $ab>1/4$); it is more elementary and seems to be forgotten.  We should also mention that Hardy proved a more general result, assuming that $|f(x)|=O(|x|^me^{-a|x|^2})$ and  $|\widehat{f}(\xi)|=O(|\xi|^me^{-b|\xi|^2})$ as $x,\xi\to\pm\infty$, he showed that $f$ is a polynomial times $e^{-a|x|^2}$.
  
  There was a search for a real variable proof of the Hardy uncertainty principle. A rather elementary (real variable)  argument, given by Terence Tao in his book \cite[\S2.6]{T1},   implies that $f$ is zero if in the statement above $ab>C_0$ for some large constant $C_0$. Another real variable proof for the case $ab>1$ is given by E. Pauwels and M. de Gosson in \cite{PdG}, surprisingly their proof  employs prolate   spheroidal  wave functions, which, in the context of time frequency analysis, first appeared in the celebrated series of works of H. Landau, H. Pollak and D. Slepian in the beginning of 1960s. The  first complete real proof for the sharp result is given  in \cite{CEKPV}.
  
  Before we exhibit the main topic of this note, the dynamical interpretation of the Hardy uncertainty principle, and give a new proof of the result,  we comment briefly on classical approaches and generalizations. 
  
  Hardy proved the theorem for the case $a=b=1/2$, which implies the general result by a simple rescaling. 
  Gilbert W. Morgan gave the following generalization of Hardy's result  already in 1934, \cite{M}. 
 
  \begin{theorem}\label{th:mor}
  Let $1<p\le 2$ and $1/p+1/q=1$, suppose that $f\in L^1(\R)$ and $|f(x)|\le Ce^{-a^p|x|^p/p}$ and $|\widehat{f}(\xi)|\le Ce^{-b^q|\xi|^q/q}$ and $ab>|\cos(p\pi/2)|^{1/p}$, then $f=0$. 
  \end{theorem}
For an interesting discussion of the Morgan theorem, extensions to functions that decay only along half-axes,
and  some remarkable related results, we refer the reader to \cite{N} and \cite{HJ}.
  
  The assumptions of  both theorems formulated above  are  point-wise bounds for a function and its Fourier transform.  In the 1980s M. Cowling and J. F. Price \cite{CP} obtained versions where the bounds  are replaced by an integral condition, the simplest version is the so-called $L^2$-Hardy uncertainty principle:
  \[e^{a|x|^2}f(x)\in L^2(\R),\quad {\text{and}}\quad e^{b|\xi|^2}\widehat{f}(\xi)\in L^2(\R)\]
 implies $f=0$ when $ab\ge 1/4$.

 Hardy's theorem can be generalized to higher dimension, the statement is exactly the same for $f\in L^2(\R^d)$. This can be deduced from the one dimensional result using the Radon transform, see \cite{SST}.
 Note that we discuss only the simplest generalization of the Hardy uncertainty principle to $\R^d$. The appealing  problem of natural higher dimensional statements is studied in \cite{BDJ, BD, D, dG}. 
  
  An interesting interpretation of  Hardy's uncertainty principle was given in the beginning of the current century, \cite{Ch, EKPV1}. It turns out that Theorem 1 is equivalent to the following statement.
 \begin{theorem}\label{th:un}
 	Let $u(t,x)$ be a solution to the free Schr\"odinger equation
 	\[\partial_t u=i\Delta u(t,x).\]
Suppose that $u\in C^1([0,T],W^{2,2}(\R^d))$  satisfies the following decay conditions \[|u(0,x)|\le Ce^{-\alpha|x|^2}\ {\text{and}}\  |u(T,x)|\le Ce^{-\beta|x|^2},\] where $\alpha,\beta>0$. \\
(i) If 	$\alpha\beta>(16T^2)^{-1}$ then $u(t,x)=0$,\\
(ii)  if $\alpha\beta=(16T^2)^{-1}$  then $u(t,x)=ce^{-(\alpha+i/(4T))|x|^2}$.
 \end{theorem} 
A real-variable proof of this theorem is due  to M.~Cowling, L.~Escauriaza, C.~E.~Kenig, G.~Ponce, and L.~Vega, \cite{CEKPV}.

In this note we first show that the uniqueness result is equivalent to Hardy's theorem  and  give a simple proof of  Theorem \ref{th:un}. The proof  involves  holomorphic functions, however the proof of part (i)  is based only on the Liouville theorem, which says that a bounded entire function is constant, the argument reminds the second proof of Theorem \ref{th:H}, given  by Hardy in \cite{H}. The proof of part (ii) requires some analysis of a singular point of a holomorphic function.  We then sketch the second proof of Hardy's theorem and give a relatively short and elementary proof of another uncertainty principle due to Beurling. The latter proof is inspired by the work of Hedenmalm, \cite{He}. To finish, we present an  overview of the recent generalizations of Theorem \ref{th:un}, which are called the dynamical versions of Hardy's uncertainty principle.

\section{Free Schr\"odinger equation}
\subsection{Solution by the Fourier transform}\label{ss:F}  In this section we present the classical formula for the solution of the Schr\"odinger equation, we provide the details for the convenience of the reader. A generalization of the result is used later in the note.
We consider the free Schr\"o\-din\-ger equation
	\begin{equation}\label{eq:S}
	\partial_t u(t,x)=i\Delta_x u(t,x),
	\end{equation}
	where $\Delta_x=\frac{\partial^2}{\partial x_1^2}
+...+\frac{\partial^2}{\partial x_d^2}$ is the Laplace operator.	It is one of  the simplest examples of a constant coefficient linear dispersive equation.
	Dispersive equations are called so since parts of solutions with different frequencies disperse with different speeds, spreading spatially. A plane wave is a solution to \eqref{eq:S} of the form
\[u_{\xi_0}(t,\xi)=\exp(ix\cdot\xi_0-it|\xi_0|^2).\]
Clearly, any superposition of the plane waves is also a solution. The plane waves satisfy $|u(t,x)|=1$. Below we  analyze solutions that decay in $x$. More precisely, we assume that $u\in C^1([0,T], W^{2,2}(\R^d))$. This smoothness assumption can be weakened but we prefer to avoid the technical details in this note.

An effective method to solve linear constant coefficient dispersive equations is by applying the Fourier transform in spatial variables. Let $\widehat{u}(t,\xi)=\mathcal{F}_xu(t,x)$, 
then \eqref{eq:S} reads
\[\partial_t\widehat{u}(t,\xi)=-i|\xi|^2\widehat{u}(t,\xi).\]
Thus the solutuon to \eqref{eq:S} with initial data $u(0,x)=u_0(x)\in L^2(\R^d)$ satisfies
\begin{equation}\label{eq:F}
	\widehat{u}(t,\xi)=e^{-it|\xi|^2}\widehat{u}_0(\xi).
\end{equation}
Hence, by the Fourier inversion formula,
\begin{align*}\label{eq:sol}
u(t,x)&=\frac{1}{(2\pi)^{d/2}}\int_{\R^d}e^{-it|\xi|^2+ix\cdot\xi}\widehat{u}_0(\xi)d\xi\\
&=\frac{1}{(2\pi)^{d}}\int_{\R^d}\int_{\R^d} e^{i(-t|\xi|^2+(x-y)\cdot \xi)}u_0(y)dyd\xi.
\end{align*}

The formula for $u(t,x)$ above  can be written as the convolution
\[u(t,x)=\int_{\R^d}u_0(y)K_t(x-y)dy,\]
where $K_t$ is the (distributional) inverse Fourier transform of the function $e^{-it|\xi|^2}$. Formally, we write
\[K_t(x)=\frac{1}{(2\pi)^{d}}\int_{\R^d}e^{i(t|\xi|^2+x\cdot \xi)}d\xi,\]
although the integral does not converge. To make sense of the integral, let
\[K^{\varepsilon}_t(x)=\frac{1}{(2\pi)^{d}}\int_{\R^d}e^{i(t|\xi|^2+x\cdot \xi)}e^{-\varepsilon|\xi|^2}d\xi.\]
Then it is easy to see that
	\[K^\varepsilon_t(x)=
\frac{1}{(4\pi(\varepsilon+it))^{d/2}}e^{-|x|^2/(4(\varepsilon+it))}.\]
	The limit of $K^\varepsilon_t(x)$ as $\varepsilon\to 0$ exists and is equal to
	\[K_t(x)=\frac{1}{(4\pi it)^{d/2}}e^{-|x|^2/(4it)}.\]
	Therefore the solution to the Schr\"odinger equation is given by
	\begin{equation}\label{eq:Ss}
	u(t,x)=\frac{1}{(4\pi it)^{d/2}}\int_{\R^d}e^{i|x-y|^2/(4t)}u_0(y)dy.
	\end{equation}
	We note that if $k_t$ denotes the standard heat kernel, then formally $K_t=k_{it}$.
	
	\subsection{Uniqueness for the free Schr\"odinger evolution and Hardy's theorem} Using the integral formula for the solution \eqref{eq:Ss}, it is not difficult to see that Theorem \ref{th:H} is equivalent to Theorem \ref{th:un} with $d=1$. We show one implication,  the Hardy uncertainty principle follows from the uniqueness result for the Schr\"odinger equation. 

	Assume that Theorem \ref{th:un} is true and let $f$ be a function as in the Hardy theorem. We define
	\[u(t,x)=\frac{1}{(4\pi it)^{1/2}}\int_{\R}e^{i|x-y|^2/(4t)-i|y|^2/4}f(y)dy,\]
	for $t>0$. Since $f$ is decaying fast the function $u(t,x)$ is smooth. Then,  differentiating the integrand, we see that $\partial_tu=i\Delta_x u$. Moreover, by taking the limit as $t\to 0$, we get $u(0,x)=e^{-i|x|^2/4}f(x)$.  Furthermore, 
	\[u(1,x)=\frac{e^{i|x|^2/4}}{(4\pi i)^{1/2}}\widehat{f}(x/2).\]
	The assumptions in the Hardy theorem can  now be translated to
	\[|u(0,x)|\le Ce^{-a|x|^2},\quad |u(1,x)|\le Ce^{-b|x|^2/4}.\]
	Now applying Theorem \ref{th:un} with $T=1$ we conclude the argument.
	
	The reverse implication can be shown in a similar way.
	
\subsection{A proof of the uniqueness theorem}
We now give  a relatively elementary proof of Theorem \ref{th:un}.  The main  idea is to consider the family of partial differential equations $\partial_t u=z\Delta_x u$ with complex parameter $z$. When $z=\pm 1$ we get the heat and the backward heat equations, while $z=i$ corresponds to the Schr\"odinger equation. Computations, similar to ones  presented in Section \ref{ss:F}, show that the fundamental solution is  
	\[k_t(z)(x)=(4\pi z t)^{-d/2} e^{-|x|^2/(4zt)}.\]
	Thus  for a fast decaying initial condition $u_0(x)$ the solution to the equation is given by $u(t,x)=u_0*k_t(z)$, so $k_t(z)=: k_{tz}$ is a complex  extension of the heat kernel. 
		
	Assume now that 
	\[|u_0(x)|=|u(0,x)|\le e^{-\alpha|x|^2}.\]
		We start with the initial condition $u(0,x)=u_0(x)$ that decays fast and we solve the generalized heat equation. We see that the heat equation itself is solvable (it corresponds to $z$ real and positive) as is the Schr\"odinger  equation (corresponding to pure imaginary $z$), but the backward heat equation cannot be solved in general, and our function is not defined for small real negative $z$. 
We consider the function
	\[F(z,x)=\frac{1}{(4\pi z)^{d/2}}\int_{\R^d}e^{-|x-y|^2/(4z)}u_0(y)dy=k_z\ast u_0,\]
	for $z\in \Omega_0=\{z\in\mathbb{C}: \Re(-1/(4z))-\alpha<0\}$. Solving the last inequality for $z$, we see that the integral above converges uniformly on compact subsets of the domain 
	 \[
	\Omega_0=\{z\in\C: |z+1/(8\alpha)|>1/8\alpha\}.\]
	The function  $F^2(z,x)$ is a holomorphic function of $z$ in $\Omega_0$,  when $x\in\R^d$ is fixed. Note that we take the square of $F$ to avoid the branching of $\sqrt{z}$.

Now, we start with $u(T,x)=u_1(x)$ and define
	\begin{equation*}
	G(z,x)=	k_{z-iT}\ast u_1=\frac{1}{(4\pi (z-iT))^{d/2}}\int_{\R^d}e^{-|x-y|^2/(4(z-iT))}u_1(y)dy.
	\end{equation*}
	Using the decay of $u_1$ we see that $G^2(z,x)$ is well defined and holomorphic in the domain
	\[\Omega_1=\{z\in \C:|z-iT+1/8\beta|>1/(8\beta)\}.\]
	Moreover $G(it,x)=u(t,x)$ when $t\in(0,T)$. Hence the holomorphic functions $F^2(\cdot,x)$ and $G^2(\cdot, x)$ coincide on the interval  $(0,T)$.  Therefore $F^2(\cdot, x)$ is extended to a holomorphic function on $\Omega_0\cup\Omega_1$. 	
	
	To simplify the notation, we denote $(8\alpha)^{-1}=A$ and $(8\beta)^{-1}=B$. Then the  complements  of $\Omega_0$ and $\Omega_1$  are circles with the radii $A$ and $B$, while the distance between the centers is $\sqrt{T^2+(A-B)^2}$. 
	
	If $AB<T^2/4$ (which is equivalent to $16\alpha\beta>T^{-2}$)  then the circles do not intersect. Thus $F^2(z,x)$ extends to an entire function in $z$ for each fixed $x$. It also satisfies
	\begin{multline}\label{eq:Fest}
	|F^2(z,x)|\le \frac{C}{(4\pi|z|)^{d}}\left(\int_{\R^d} e^{-\Re(|x-y|^2/(4z))}e^{-\alpha|y|^2}dy\right)^2\\= \frac{C}{(4|z|(\alpha+\gamma))^{d}}e^{-2\gamma\alpha|x|^2/(\gamma+\alpha)},
	\end{multline}
	where $\gamma=\Re(1/(4z))$. We fix $x$ and note that  $F^2(z,x)$ is uniformly  bounded as $|z|>1/\alpha$. Then,  by the Liouville theorem,  $F^2(z,x)$ is a constant function in $z$ for each $x$. This means that $\partial_t u=0$ and thus $\Delta u=0$. There are no non-zero decaying harmonic functions, therefore $u(t,x)=0$.  
	
This proof of part (i) uses only the facts that the function $e^{cz}$ satisfies the mean value property and that a bounded function satisfying the mean value property on the whole plane is a constant. 
An elementary proof of the latter can be found in \cite{Ne}.

	\def\xaxis{(-5.0,0.0) -- (3.0,0.0)}
	\def\yaxis{(0,-1.0)--(0,3.8)}
	\def\Sa{(-1,0) circle (1cm)}
	\def\Sb{(-1.5635, 2.5) circle (1.5635cm)}
	
	\begin{figure}
	
\begin{center}	
	\begin{tikzpicture}[remember picture]
		\draw[thick,->] \xaxis;
		\draw[thick,->] \yaxis;
		\draw[blue, thick, name path = S1]\Sa;
		\draw[green, thick, name path = S2] \Sb;
	\path [name intersections={of=S1 and S2, by=Z}];
		\node[label={$z_0$}] at (Z) {$\bullet$};
		\node[label=2:{$iT$}] at (0,2.5) {$\bullet$};
		\node at (-1.5635, 2.5) {$\cdot$};
		\path [name path=aux] (Z) circle [radius=0.3bp];
		\draw [red, thick, name intersections={of=S1 and aux}] (Z) -- ($(intersection-1)!3cm!(intersection-2)$);
			\draw [red, thick, name intersections={of=S2 and aux}] (Z) -- ($(intersection-2)!3cm!(intersection-1)$);
			\node[label={$z$ plane}] at (2,3) {};
 	\end{tikzpicture}
	
	\caption{Tangent circles $\partial\Omega_0$ and $\partial\Omega_1$ and their common tangent line $l$ for the case $AB=T^2/4$,  $z$-plane}
	\label{f:1}
	\end{center}
	\end{figure}
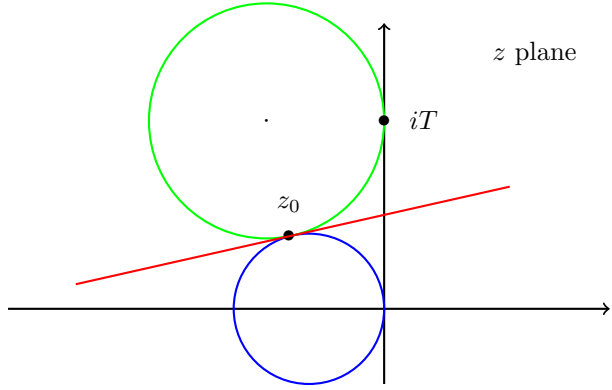

Now assume that   $16\alpha\beta=T^{-2}$, i.e., $AB=T^2/4$, then  the circles $\partial\Omega_0$ and $\partial\Omega_1$ touch at one point, which we denote by  $z_0$, see Figure \ref{f:1}.
Thus $F^2(z,x)$ is a holomorphic function in $\C\setminus\{z_0\}$.  We consider $x=0$ and claim that $F^2(z,0)$ has a pole  at $z_0$.  
	To prove that, we draw the common tangent line $l$ to the circles $\partial\Omega_0$ and $\partial\Omega_1$, and consider the images of this line under the transformations $\zeta=z^{-1}$ and $\eta=(z-iT)^{-1}$. These are circles $\omega_0$ and $\omega_1$ passing through the origin, while the images of the circles $\partial\Omega_0$ and $\partial\Omega_1$ under those two respective transformations are vertical lines $l_0$ and $l_1$ tangent to $\omega_0$ and $\omega_1$, see Figure \ref{f:2}. We see that $\omega_0$ is defined by the equation 
	\[\Re(\zeta-\zeta_0)=|\zeta-\zeta_0|^2/(2r_0),\]
	where $\zeta_0=z_0^{-1}$ and $r_0$ is the radius of $\omega_0$.  Let $z$ be a point close to $z_0$ lying above the line $l$ (on the other side of the line $l$ than $\partial\Omega_0$). Then $\zeta=z^{-1}$ lies inside the disk bounded by $\omega_0$ and we have the following inequality
	\begin{equation}\label{eq:ReC}
	\Re(\zeta-\zeta_0)\ge c|\zeta-\zeta_0|^2\ge c_1|z-z_0|^2,
	\end{equation}
	where $c=(2r_0)^{-1}$  and $c_1=c|z_0|^{-4}/2$. The estimate \eqref{eq:Fest} implies
	\[|F^2(z,0)|\le C|z-z_0|^{-2d}\]
	when $z$ is in the half-plane above the line $l$. For the other half-plane we repeat the argument, using the  function $G^2$, and conclude that $F^2(z,0)$ has a pole at $z_0$ of order less than or equal to $2d$. 
	
	Similarly, we consider the  functions
	\begin{equation*}
	F_j(z,x)=\partial F(x,z)/\partial x_j=\frac{1}{2z(4\pi z)^{d/2}}\int_{\R^d} e^{-|x-y|^2/(4z)}(y_j-x_j)u_0(y)dy,\quad j=1,...,d.
	\end{equation*}
	Then  each $F_j^2(z,x)$ extends to a holomorphic function in $\C\setminus\{z_0\}$ and $F_j(z,0)$ has a pole at $z_0$. An estimate of $F_j(z,0)$  gives
	\[|F_j^2(z,0)|\le C|z-z_0|^{-2d-1}.\]
	
	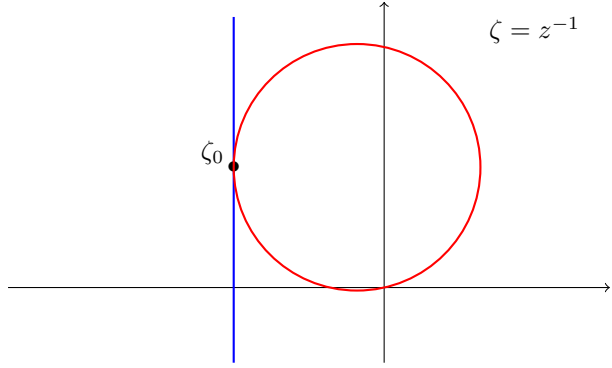
\begin{figure}\label{f:2}
	\begin{center}
		\begin{tikzpicture}
	\draw[->] \xaxis;
	\draw[->] \yaxis;
	\path [name path =S1]\Sa;
	\path[name path = S2]\Sb;
	\path [name intersections={of=S1 and S2, by=Z}];
	\coordinate (o) at (0,0);
	\coordinate (z) at (Z);
	\coordinate(w) at (-2,1.6);
	\coordinate (v) at (-0.72,1.6);
	\draw[thick, blue, name path = L1] (-2,-1)--(-2,3.6);
	\path [name path=L] (0,0) -- ($(o)!5cm!(z)$);
	\path [name intersections={of=L1 and L, by=W}];
	\node[label={[left]$\zeta_0$}] at (w) {$\bullet$};
	\draw[red, thick] (-0.36,1.6) circle (1.64cm);
	\node[label={$\zeta=z^{-1}$}] at (2,3) {};
	\end{tikzpicture}
	\caption{Circle $\omega_0$ and  tangent line $l_0$  in $\zeta=1/z$-plane}
	\label{f:2}
	\end{center}
	\end{figure}
	
	Finally,  consider  $\zeta=\zeta_0+t$, where $t>0$ is real and small. For this case the inequality \eqref{eq:ReC} can be replaced by $\Re(\zeta-\zeta_0)=|\zeta-\zeta_0|$.  Then, repeating the argument above and taking   $z=1/\zeta$, we see that $(z-z_0)^dF^2(z,0)$ and $(z-z_0)^{d+1}F_j^2(z,0)$ are bounded along the curve $z=z_0(1+tz_0)^{-1}$, $t>0$. Thus $F^2(z,0)$ has a  pole at $z_0$ of order not exceeding $d$, while for each $F_j^2(z,0), j=1,...,d$, the order of this  pole  does not exceed $d+1$. 
	
	We assume first that $d=1$. 	To finish the proof of the end-point case we use the Hermite functions,
		\[\psi_n(t)=e^{t^2/2}\frac{\partial^n}{\partial t^n}e^{-t^2}=H_n(t)e^{-t^2/2},\]
	which form an orthogonal basis for $L^2(\R)$.		
More generally, for any complex number $\gamma$ with $\Re\gamma>0$ we may define the generalized Hermite functions 
\[\psi^{(\gamma)}_n(t)=H_n(\sqrt{2\gamma} t)e^{-\gamma t^2},\]
which still form an orthogonal basis for $L^2(\R)$.

First we consider $F^2(z,0)$.  This is a holomorphic function in $\C\setminus\{z_0\}$ that tends to  zero at infinity and has a simple pole at $z_0$, thus
	\[F^2(z,0)=b(z-z_0)^{-1}.\]
	Hence
	\[\int_{-\infty}^\infty e^{-y^2/(4z)}u_0(y)dy=cz^{1/2}(z-z_0)^{-1/2}.\]
	A simple computation shows that
\begin{equation}\label{eq:gH}
	\int_{-\infty}^\infty e^{-y^2/(4z)}e^{-\gamma y^2}dy=2\sqrt{\pi}\frac{z^{1/2}}{(1+4z\gamma)^{1/2}}.
	\end{equation}
We choose  $\gamma=-1/(4z_0)=-\zeta_0/4$ and see that for some constant $c_0$ and every $\zeta$
\[c_0\int_{-\infty}^\infty e^{-y^2\zeta/4}e^{-\gamma y^2}dy=\int_{-\infty}^\infty e^{-y^2\zeta/4} u_0(y)dy.\]
This means that all even moments of $u_0$ are equal to the corresponding moments of $c_0e^{-\gamma |y|^2}$ and thus $u_0(y)+u_0(-y)=2c_0e^{-\gamma y^2}$. 

Then, similarly, we consider $F^2_1(z,0)$. We have
\[F_1^2(z,0)=b_2(z-z_0)^{-2}+b_1(z-z_0)^{-1}.\]
 On the other hand
\begin{equation}\label{eq:gH1}
\int_{0}^\infty y^ke^{-y^2/(4z)}e^{-\gamma y^2}dy=c_k\frac{z^{(k+1)/2}}{(1+4z\gamma)^{(k+1)/2}}.
\end{equation}
Representing $yu_0(y)$ as the series in $\psi^{(\gamma)}_n$, we 
conclude that 
\[yu_0(y)-yu_0(-y)=2(c_1+c_2y) e^{-\gamma y^2}.\]
 Now, taking $y\to 0$ and using that $u_0(y)=G(0,y)$ is a continuous function, we see that  $c_1=c_2=0$. Thus $u_0$ is even and $u_0(y)=c_0e^{-\gamma |y|^2}$.  It is not difficult to check  that $\gamma=\alpha+i/4T$. This concludes the proof of Theorem \ref{th:un} for the case $d=1$.

To complete the proof in higher dimensions we consider $F(z,x)$ and all its partial derivatives in the spatial variables  at $x=0$. Rewriting the integral in polar coordinates, we have
\[F(z,0)=\frac{1}{(4\pi z)^{d/2}}\int_0^{\infty} r^{d-1}\int_{S^{d-1}} u_0(ry')d\sigma(y') e^{-r^2/(4z)}dr.\]
Let  $\Phi(r)=r^{d-1}\int_{S^{d-1}}u_0(ry')d\sigma(y')$. The identity \eqref{eq:gH1} and the fact that $F^2$ has a pole at $z_0=-1/(4\gamma)$ of order not exceeding $d$ imply that  
\begin{equation*}
\Phi(r)=\sum_{l=0}^{d-1}c_lr^{l}e^{-\gamma r^2}.
\end{equation*}  
Moreover, since $\Phi(r)$ has a zero of order $d-1$ at zero, we conclude that $\Phi(r)=cr^{d-1}e^{-\gamma r^2}$. On the other hand, looking at the partial derivatives of $F$ we see that for any  homogeneous polynomial $p(y)$ of degree $k$, 
\[\Phi_{p}(r)=r^{d-1+k}\int_{S^{d-1}} p(y') u_0(ry')d\sigma(y')\] 
is  a linear combination of the form $\sum_0^{d-1+k} c_lr^{l}e^{-\gamma r^2}$. If $\int_{S^{d-1}}p(y')d\sigma(y')=0$ then $\Phi_p(r)=0$ since its zero at the origin is of order larger than $d-1+k$. Therefore $u_0$ is orthogonal to all polynomials with zero mean on each sphere centered at the origin. This implies that $u_0$ is a constant on each such sphere and thus $u_0(y)=ce^{-\gamma|y|^2}$.

\subsection{Heat equation}
We saw that the Schr\"odinger equation and the heat equation are close relatives. Therefore, it is natural that the Hardy uncertainty principle implies a uniqueness result for the heat equation. 
\begin{theorem}\label{th:heat}
 Let $u(t,x)\in C^1([0,T], W^{2,2}(\R^d))$ be a solution to the heat equation $\partial_t u=\Delta_x u$. Suppose  that $u(0,x)\in L^1(\R^d)$ and
$|u(T,x)|\le e^{-\delta |x|^2}.$
If $\delta\ge 1/(4T)$ then $u=0$.
\end{theorem}
The case $\delta=1/4T$ corresponds to the situation $u(0,x)$ is the Dirac delta function. The fact that  the Hardy uncertainty principle implies Theorem \ref{th:heat} follows by applying the Fourier transform in variable $x$, which gives
\begin{equation}\label{eq:heat1}
\widehat{u}(t,\xi)=e^{-t|\xi|^2}\widehat{u}(0,\xi).
\end{equation}
Thus, if the initial data $u_0(x)=u(0,x)\in L^1(\R^d)$ then $|\widehat{u}(T,\xi)|\le Ce^{-T|\xi|^2}$, combined with the decay condition for $u(T,x)$, it implies  that $u(T,x)=0$ if $\delta>1/(4T)$ and $u(T,x)=c_0e^{-\delta|x|^2}$ if $\delta=1/(4T)$. The latter implies $\widehat{u}_0(\xi)=c$ and $u$ is a multiple of the Dirac  delta function.

We can also prove Theorem \ref{th:heat} using the approach suggested in the  previous section. The condition $|u(T,x)|\le e^{-\delta|x|^2}$ implies that the function
\[\tilde{G}^2(z,x)=(k_{z-T}*u(T,x))^2\]
is holomorphic in the domain
\[\tilde{\Omega}=\{z\in\mathbb{C}:|z-T+(8\delta)^{-1}|>(8\delta)^{-1}\}.\]
While the condition $u(0,x)\in L^1$ implies that the function 
\[\tilde{F}^2(z,x)=(k_z\ast u(0,x))^2\]
is holomorphic when $\Re(z)>0$. Moreover we know that $\tilde{F}^2(t,x)=\tilde{G}^2(t,x)$ when $t\in(0,T)$. If $\delta>1/(4T)$, the two domains cover the whole complex plane and we obtain a bounded entire function. It leads to a contradiction in the same way as above for the Schr\"odinger equation. If $\delta=1/(4T)$ then the resulting function is holomorphic in $\C\setminus\{0\}$, but  the singularity at $0$ is removable for almost every $x$ since
\[\lim_{z\to 0} \tilde{F}^2(z,x)=u^2(0,x)\]
almost everywhere. And we get a contradiction again.

We also note that Theorem \ref{th:heat} does not imply the limit case ( $ab=1/4$) in the Hardy uncertainty principle. The reason is that in general a bounded function is not a Fourier transform of an $L^1$-function. To obtain an equivalent statement, one should extend the notion of solutions of the heat equation to the case when the initial data is a measure.

\section{The second proof of Hardy and Beurling's uncertainty principle}
\subsection{On forgotten proof of Hardy} We were not able to find the second proof of Hardy or its variations in  any textbook, so we give a sketch of this proof here, as pointed out in the introduction, for the case $a=b=\frac12$. First, Hardy notes that the decay conditions on $f$ and $\widehat{f}$ imply the decay conditions on $f_e(x)=(f(x)+f(-x))/2$ and $f_o=(f(x)-f(-x))/2$ and their Fourier transforms. Next, the functions $f_1=(f_e+\widehat{f_e})/2$, $f_2=(f_e-\widehat{f}_e)/2$, $f_3=(f_o+i\widehat{f_o})/2$, and $f_4=(f_0-i\widehat{f_o})/2$ also satisfy the decay condition together with the Fourier transforms.   So one may assume that $\widehat{f}=i^kf$. 

	Let  first $f$ be even, so that $\widehat{f}=\pm f$. Hardy considers the function
	\[\lambda_f(s)=\int_0^\infty e^{-sx^2/2} f(x)\,dx,\]
	where $f$ decays as the Gaussian. Then $\lambda_f$ is a holomorphic function when $\Re(s)>-1$ and the equation $\widehat{f}=\pm f$  translates into the identity
	\[\lambda_f(s)=s^{-1/2}\lambda_f(1/s),\]
	we skip the details of choosing the right branch of the root function here. 
	
	Then the  function  $\mu(s)=\sqrt{s+1}\,\lambda_f(s)$ satisfies $\mu(s)=\mu(1/s)$ and it  can be extended to a holomorphic function in $\C\setminus\{-1\}$. Moreover, $\mu$ has a pole at $s_0=-1$. Finally, Hardy refers to the injectivity of the transform, i.e., $\lambda_f=\lambda_g$ if and only if $f=g$, and the identity for the Hermite functions
	\begin{equation*} 
		\int_{0}^\infty \psi_{2n}(t)e^{-st^2/2}dt=c_n\frac{(s-1)^{n}}{(s+1)^{n+1/2}}.
	\end{equation*}

The case $f$ is odd is not written down in \cite{H}. For this case we suggest to consider the function
\[\tilde{\lambda}_f(s)=\int_0^\infty xe^{-sx^2/2} f(x)\,dx=\pm \sqrt\frac{2}{\pi}\int_0^\infty xe^{-sx^2/2}\int_0^\infty f(y)\sin xy\, dy\,dx,\]
the second identity follows from the fact $f=\pm i\widehat{f}$. Then $\tilde{\lambda}_f(s)=s^{-3/2}\tilde{\lambda}_f(1/s)$. As before, we consider $\mu(s)=\sqrt{(s+1)}\,\tilde{\lambda}_f(s)$ that satisfies $\mu(s)=s^{-1}\mu(1/s)$. This function extends to a holomorphic function in $\C\setminus\{-1\}$ such that $|\mu(s)|\to0$ when $|s|\to\infty$. Further, $\mu$ has a pole at $s_0=-1$, and one concludes the argument by the same techniques of the even case.

\subsection{Beurling's uncertainty principle} The following version of the uncertainty principle  is due to Arne Beurling
\begin{theorem} Suppose that $f\in L^2(\R)$ and 
\[\int_\R\int_\R e^{|x\xi|}|f(x)||\widehat{f}(\xi)|\,dx\,d\xi<\infty.\]
Then $f=0$.
\end{theorem}

The theorem appeared in the collected works of Beurling, \cite{B} and dates back to the 1960s. The original proof of Beurling uses the Phragm\'en-Lindel\"of theorem and it  can be found in \cite{Ho}. Higher dimensional versions of the Beurling theorem were obtained in \cite{BDJ}. In 2012 H\aa kan Hedenmalm gave another proof and generalized the statement in \cite{He}. His result was further extended in \cite{G}. We follow the ideas  in \cite{He} to give a relatively short proof of the original statement of Beurling. Clearly, the Beurling theorem implies the $L^2$-version of the Hardy uniqueness result.

First, by taking the real and imaginary parts of $f$ we may reduce the problem to the case when $f$ is real-valued. Now, following the idea of Hedenmalm, consider the function
\[F(s)=\int_\R\int_\R e^{is x\xi}f(x)\widehat{f}(\xi)\, dx\, d\xi.\]
Then $F$ is well-defined and holomorphic in the strip $S=\{s\in\C: |\Im(s)|<1\}$. Moreover, by the monotone convergence theorem, $F$ is continuous on $\overline{S}$. For real $s$, we have
\[F(s)=\sqrt{\pi/2}\int_\R f(x)f(s x)\, dx,\]
we have used that $f,\widehat{f}\in L^1(\R)$.
Then $F(s)=s^{-1}F(1/s)$ for $s\in\R\setminus\{0\}$.   
We obtain that $F$ can be extended to a holomorphic function on $\C\setminus\{\pm i\}$. The singularities at $s=\pm i$ are removable since the function is continuous at these points. Finally, the functional equation $F(s)=s^{-1}F(1/s)$ and the fact that $F$ is bounded near the origin imply that $|F(s)|\to 0$ when $|s|\to\infty$. Thus $F=0$. In particular,
\[F(1)=\sqrt{\frac{\pi}{2}}\int_\R f^2(x)\, dx=0.\]
Finally, since $f$ is real-valued, we conclude that $f=0$.

\section{Recent versions of the uniqueness theorem}
We now return to the dynamical versions of the uncertainty principles. In the last 15 years the uniqueness results for the free Schr\"odinger and  heat equations were generalized to a large class of evolutions. We give an overview of some of these results in this section.
\subsection{Schr\"odinger and heat equations with a potential}\label{ss:mS}
First, we consider the Schr\"odinger equation with a potential,
\begin{equation}\label{eq:sch}
\partial_t u(t,x)=i(\Delta u+Vu).
\end{equation}
In a series of articles, Luis Escauriaza, Carlos E. Kenig, Gustavo Ponce, and Luis Vega, \cite{EKPV1, EKPV3, EKPV4, EKPV, EKPV2}, generalized the uniqueness result for the case when $V$ is a bounded potential satisfying one of the following conditions:\\
(i) $\lim_{R\to\infty}\int_0^T\sup_{|x|>R}|V(t,x)|dt=0$,\\
(ii) $V(t,x)=V_1(x)+V_2(t,x)$, where $V_1$ is real-valued (and does not depend on $t$) and $V_2$ satisfie, for some positive $\alpha$ and $\beta$,
\[
\sup_{[0,T]}\| e^{\alpha\beta T^2|x|^2/(\sqrt{\alpha}t+\sqrt{\beta}(T-t))^2}V_2(t)\|_{L^\infty(\mathbb{R}^n)}<+\infty.
\]
\begin{theorem}\label{th:m}
	Let $u\in C([0,T], L^2(\R^d))$ be a solution to \eqref{eq:sch}, where $V$ satisfies either (i) or (ii). If $|u(0,x)|\le Ce^{-\alpha|x|^2}$ and $|u(1,x)|\le Ce^{-\beta|x|^2}$ with $\alpha\beta>1/(16T^2)$ then $u=0$.
\end{theorem}

Note that the condition on $\alpha\beta$ is sharp! The result is further generalized to semi-linear equations and covariant Schr\"odinger evolution in \cite{EKPV} and \cite{BFGRV, CF}, and to Navier-Stokes equation in \cite{DHS}.

We outline the proof of Theorem \ref{th:m}. First it suffices to consider the case when $\alpha=\beta$, the Appell transform reduces the general case to this one.  We renormalize the solution and assume that $T=1$. The first step is to show logarithmic  convexity of some weighted norm of the solution, the method can be compared to the one  used by Shmuel Agmon for elliptic equations in 1960s, see \cite{A}. For each $t\in[0,1]$ and $\xi\in S^{d-1}$ we define 
\[H(t)=\int_{\R^d}| e^{\mu|x+Rb(t)\xi|^2}u(t,x)|^2dx,\]
where $b(t)=16\mu t(1-t)$. The derivative of $v(t,x)=e^{\mu|x+Rb(t)\xi|^2}u(t,x)$ in $t$  is written as the sum of a symmetric and anti-symmetric operator,
\[\partial_t v=(S+A)v.\]
Then a straightforward calculation implies that 
\[(\log H(t))''\ge 2\langle (SA-AS)v,v\rangle.\]
Careful estimates on $SA-AS$ show that  $(\log H(t))''\ge -16\mu R^2-C_V$, where $C_V$ denotes a constant that depends on the potential. Therefore
\begin{equation}\label{eq:conv}
H(t)\exp(-32\mu R^2t(1-t))\le C_VH(0)^{1-t}H(1)^{t}.
\end{equation}
The right hand side does not depend on $R$, while in the left hand side for $t=1/2$ the weight (with which $u^2$ is integrated) is 
\[\exp(2\mu|x+4\mu R\xi|^2-8\mu R^2).\]
We look at the coefficient in front of $R^2$, if  $32\mu^3>8\mu$ it is positive and thus we see that $u(1/2, x)=0$ for almost each $x$, by letting $R\to \infty$. Then $u\equiv 0$. This formal computation can be justified if $H(0)$ and $H(1)$ are finite. This  proves Theorem \ref{th:m} when $\alpha=\beta>1/2$.  

To extend the result for the range $\alpha=\beta>1/4$, Escauriaza, Kenig, Ponce, and Vega developed an ingenious bootstrapping argument. To  sketch their argument, we write \eqref{eq:conv} as
\[ \int_{\R^d}|u(t,x)|^2e^{2\mu|x|^2+4R \mu b(t)x\cdot\xi-2R^2b(t)(1-\mu b(t))}dx\le C_VH(0)^{1-t}H(1)^{t}.\]
Under the assumption $\alpha=\beta\le1/2$ a formal integration of the last inequality with respect to $R$ leads to
\[\int_{\R^d}|u(t,x)|^2e^{2a_1(t)|x|^2}dx\le C_VH(0)^{1-t}H(1)^t,\]
for $a_1(t)=\mu/(1-\mu b(t))$. Notice that $a_1(1/2-t)=a_1(1/2+t)$, $a_1(0)=a_1(1)=\mu$ and $a_1(t)>\mu$ when $t\in(0,1)$, which shows that the solution $u$ decays faster at $(0,1)$ than at the endpoints. Next, one can construct a positive function $b_1(t)$ such that $b_1(0)=b_1(1)=0$ and so that
\[H_1(t)=\int_{\R^d}| e^{a_1(t)|x+Rb_1(t)\xi|^2}u(t,x)|^2dx,\]
satisfies
\begin{equation}\label{eq:convk}
H_1(t)\exp(-2R^2b_1(t))\le C_VH_1(0)^{1-t}H_1(1)^{t}=C_VH(0)^{1-t}H(1)^{t}.
\end{equation}
Note that this is again  \eqref{eq:conv} but $\mu$ and $b$ are replaced by  $a_1$ and $b_1$. A similar study as before tells us that $1-a_1(1/2)b_1(1/2)\le 0$ implies $u\equiv0$, while otherwise we can integrate again to improve the decay at $(0,1)$. This self-improvement can be repeated several times, resulting in a sequence of functions
\begin{equation}\label{eq:ak}a_{k+1}(t)=\frac{a_k(t)}{1-a_k(t)b_k(t)},\ a_0(t)=\mu\end{equation}
such that
\[\mu<a_1(t)<...<a_k(t),\quad t\in(0,1). \]
On each step the new function satisfies $a_k(1/2-t)=a_k(t+1/2)$, $a_k(0)=a_k(1)=\mu$, and
\[\|e^{a_k(t)|x|^2}u(t,x)\|_2^2\le C_VH(0)^{1-t}H(1)^t.\]
As for the functions $b_k(t)$, they are constructed from $a_k(t)$ in such a way that at each step relation \eqref{eq:convk} is satisfied for the pair of functions $a_k$ and $b_k$. More precisely, as shown in \cite{EKPV4}, $b_k(t)$ is the solution to
\[
\left\{\begin{array}{l}\ddot b_k=-\frac{1}{a_k^2}\big(\ddot a_k+32a_k^3-\frac{3(\dot a_k)^2}{2a_k}\big),\\
b_k(0)=b_k(1)=0.\end{array}\right.
\]

If, for some $k$, we have $1-a_k(1/2)b_k(1/2)\le0$, which translates in a condition on parameter $\mu$, the iterative argument stops and we reach a contradiction implying $u\equiv0$. Otherwise, the process is infinite and the limit function $a(t)=\lim_{k\to\infty} a_k(t)$ exists. Since \eqref{eq:ak} implies $b_k(t)=(a_{k+1}-a_k)/(a_k a_{k+1})$, the functions $b_k$ will converge to 0 and, from the differential equation satisfied by $b_k$, one can deduce that the limit function $a(t)$ satisfies 
\[
\left\{\begin{array}{l} \ddot a+32a^3-\frac{3(\dot a)^2}{2a}=0,\\
a(0)=a(1)=\mu.\end{array}\right.
\]

Solving the ODE under the constraint $a(1/2-t)=a(1/2+t)$ leads to 
\[
a(t)=\frac{C}{4\big(1+(t-1/2)^2C^2\big)}
\]
for some $C>0$. 
Computing the maximum in $C$ of $\mu=a(0)=C/(4+C^2)$,  we see  that $\mu$ must be less than $1/4$. Then Theorem \ref{th:m} follows.

A similar strategy gives a powerful generalization of Theorem \ref{th:heat},\cite{EKPV2}.
\begin{theorem}
	Let $V(t,x)\in  L^\infty(\R\times\R^d)$ and $u$ be a solution to the equation
	\[\partial_t u=\Delta_x u+Vu,\]
	$u\in L^\infty([0,T], L^2(\R^d))\cap L^2([0,T], H^1(\R^d)])$.  If $|u(T,x)|\le e^{-\delta|x|^2}$ and $\delta>1/\sqrt{T}$, then $u=0$.
\end{theorem}

A natural question is what decay  a stationary solution to the Schr\"odinger equation may have. The question was asked by E. M. Landis in 1960 (see \cite{BK,KL}), who conjectured that if $V\in L^\infty(\R^d)$, $\Delta u+Vu=0$ in $\R^d$, and $|u(x)|\le C\exp(-|x|^{1+c})$ with $c>0$, then $u\equiv 0$. The conjecture was disproved by V. Z. Meshkov in \cite{Me}, who constructed an example of a complex valued $u$ and $V$ such that $|u(x)|\le \exp(-|x|^{4/3})$ and proved that there are no solution with a faster decay. A remaining question is if the Landis conjecture holds under the assumption that  $V$ is real valued. In spite of some recent progress \cite{LMNN},  this is an open problem in dimensions $d\ge 3$.

\subsection{Discrete evolutions} Another twist of the uniqueness results for Schr\"odinger equation was given  in \cite{FBV, JLMP, FB1, FB2}, where uniqueness theorems are obtained for  the  discrete equation. Let $\Delta_d$ be the usual discrete  Laplacian on $\Z^d$. We consider the equation
\begin{equation}\label{eq:DS}
\partial_tU(t, n)=i(\Delta_d U(t,n)+V(t,n)U(t,n)),
\end{equation}
where $n\in \Z^d$ and $V$ is a bounded potential. The uniqueness results say  that  a solution to the discrete Schr\"odinger equation which decays fast at two times is trivial. To find the optimal decay, we consider the free evolution with $V=0$.  In dimension $d=1$, there is a solution $U_0(t,n)=i^{-n}e^{-2it}J_n(1-2t)$, where $J_n$ is the Bessel function, and it has optimal decay at $t=0$ and $t=1$. The role of the Gaussian is now played by the Bessel function. This fact is related to different behavior of the heat kernels:  for the continuous case the standard heat kernel is $k(1,x)=(4\pi)^{-1/2}\exp(-x^2/4)$, while for the discrete case the heat kernel is $K(1,n)=e^{-1}|I_n(1)|\asymp e^{-1}(n!2^n)^{-1}$, where $I_n$ are the modified Bessel functions, $I_n(z)=(-i)^nJ_n(iz)$.

\begin{theorem}
Let $U(t,n)$ be a solution to  \eqref{eq:DS}, with $V\equiv 0$, on $[0,1]\times\Z $. Suppose that
\[
|U(0,n)|+|U(1,n)|\le\frac{C}{\sqrt{|n|}}\left(\frac{e}{2|n|}\right)^{|n|},\quad n\in\Z\setminus\{0\}.
\]
Then $U(t,n)=Ci^{-n} e^{-2it}J_n(1 - 2t)$. In particular, a solution to the free discrete Schr\"odinger  equation cannot decay faster than $J_n(1)$ both at $t=0$ and $t=1$. 
\end{theorem}
The idea of the proof is to consider the function $\psi(t,z)=\sum_{-\infty}^\infty U(t,n)z^n$. It is not difficult to show that it is defined on the unit circle $|z|=1$, Moreover the decay of $U(0,l)$ and $U(1,l)$ shows that $\psi(0,z)$ and $\psi(1,z)$ are entire functions.  The equation \eqref{eq:DS} implies
\[
\psi(t,z)=e^{i(z+z^{-1}-2)t}\psi(0,z),
\]
and $\psi(t,z)$ extends to an entire function for any $t\in[0,1]$.
Careful analysis of this function and application of the Phragm\'en--Lindel\"of theorem finishes the proof. It would be interesting to find a real-variable, or at least more elementary, proof.

This result was  generalized to special classes of time-independent potentials. General bounded potentials were considered in \cite{JLMP} (in dimension $d=1$) and \cite{FBV} (in arbitrary dimension).  The result is as follows.

\begin{theorem}\label{th:DV}
 Let $U(t,n) \in C^1([0,1]:\ell^2(\mathbb{Z}^d))$ be a solution to  \eqref{eq:DS} on $[0,1]\times\Z ^d$. Suppose that $\|V\|_\infty\le 1$. There exists constant $\gamma$ such that if 
\[
|U(0,n)|+|U(1,n)|\le C\exp(-\gamma |n|\log |n|),\quad n\in\Z^d\setminus\{0\}.
\]
then $U=0$.
\end{theorem}

The approach in \cite{JLMP} follows the scheme of \cite{EKPV} described in the first step of the proof of Theorem \ref{th:m} in Section \ref{ss:mS}. We describe the details of \cite{FBV}. The idea  is to make use of the following result, known in the literature as Carleman-type inequality, whose proof relies on the computation of a commutator between a symmetric and an anti-symmetric operator. In what follows $\|\cdot\|_2$ stands for $\|\cdot\|_{L^2([0,1],\ell^2(\mathbb{Z}^d))}$, and $\|\cdot\|_\infty$ will represent the supremum norm. 
\begin{lemma}
\label{lem21}
Let $\varphi:[0,1]\rightarrow\mathbb{R}$ be a smooth function and $\gamma>\frac{\sqrt{d}}{2}$. There exists $R_0=R_0(d,\|\varphi'||_\infty+\|\varphi''\|_\infty,\gamma)$ and $c=c(d,\|\varphi'\|_\infty+\|\varphi''\|_\infty)$ such that, if $R>R_0$, $\alpha\ge \gamma R\log R$ and $g\in C_0^1([0,1],\ell^2(\mathbb{Z}^d))$ has its support contained in the set
\[
\{(t,n):|n/R+\varphi(t) e_1|\ge 1\}
\]
then
\[\begin{aligned}
\sqrt{\sinh(2\alpha/R^2)}\sinh(2\alpha/\sqrt{d}R)\|e^{\alpha\left|\frac{n}{R}+\varphi(t) e_1\right|^2}g\|_2\le c \|e^{\alpha\left|\frac{n}{R}+\varphi(t) e_1\right|^2}(i\partial_t+\Delta_d)g\|_2.
\end{aligned}\]
\end{lemma}

Thanks to this inequality, one can deduce lower bounds for nontrivial solutions of \eqref{eq:DS} with a general bounded potential. In order to do that, consider the following cut-off functions
\[
\theta^R(x)=\begin{cases}1, |x|\le R-1,\\0, |x|\ge R,\end{cases}\ \  \mu(x)=\begin{cases}1,|x|\ge 2,\\0,|x|\le 1,\end{cases}\ \  \varphi(t)=\begin{cases}3,t\in[\frac38,\frac58],\\0,t\in [0,\frac14]\cup[\frac34,1],\end{cases}
\]
and define $g(t,n)=U(t,n)\theta^R(n)\mu\big(\frac{n}{R}+\varphi(t)e_1\big)$. By means of the Leibniz rule, and carefully studying the size of the weight $e^{\alpha\left|\frac{n}{R}+\varphi(t) e_1\right|^2}$ in the support of the derivatives of the cut-off functions, one can check that
\begin{equation}\begin{aligned}\label{eq:carl}
\sqrt{\sinh(2\alpha/R^2)}&\sinh(2\alpha/\sqrt{d}R)\|e^{\alpha\left|\frac{n}{R}+\varphi e_1\right|^2}g\|_2\le c\|e^{\alpha\left|\frac{n}{R}+\varphi e_1\right|^2}(i\partial_t+\Delta_d)g\|_2\\ \le& c\left(\|e^{\alpha\left|\frac{n}{R}+\varphi e_1\right|^2}g\|_2+e^{16\alpha}\lambda(R)+e^{4\alpha}\|U\|_2\right),
\end{aligned}\end{equation}
where $\lambda(R)= \left(\int_{0}^{1}\sum_{R-2\le |n|\le R+1}|U(t,n)|^2\right)^{1/2}$. The fact that $\alpha$ needs to be larger that $\gamma R\log R$ implies that for $R \ge R_0$ depending only on the dimension, the first term in the right-hand side can be absorbed in the left-hand side (one can check that the product of $\sinh$ functions increases with $R$). On the other hand, if we assume $\int_{1/2-1/8}^{1/2+1/8}|U(t,0)|^2\,dt\ge1,$ the norm in the left-hand side is bounded by
\[
\|e^{\alpha\left|\frac{n}{R}+\varphi e_1\right|^2}g\|_2\ge e^{9\alpha},
\]
since $g(t,0)=U(t,0)$ if $t\in[1/2-1/8,1/2+1/8]$, and in that the region the weight is exactly $e^{9\alpha}$. So for $R\ge R_0$ depending on $\|U\|_2$ the last term in the right-hand side of \eqref{eq:carl} can also be absorbed and we get
\[
c\left(\int_{0}^{1}\sum_{R-2\le |n|\le R+1}|U(t,n)|^2\right)^{1/2}\ge  e^{-5\alpha}= e^{-cR\log R}
\]
after choosing $\alpha$ appropriately. This proves the following lower bound.

\begin{theorem}\label{th:DSV}
Let $U\in C^1([0,1]:\ell^2(\mathbb{Z}^d))$ satisfy  \eqref{eq:DS}. Assume that
\[
\int_0^1\sum_{n\in\mathbb{Z}^d}|U(t,n)|^2\,dt\le A^2,\ 
\int_{1/2-1/8}^{1/2+1/8}|U(t,0)|^2\,dt\ge1,
\]
and
\[
\|V\|_{\infty}=\sup_{t\in[0,1],j\in\mathbb{Z}^d}\{|V(t,n)|\}\le 1,
\]
then there exists $R_0=R_0(d,A)>0$ and $c=c(d)$ such that for $R\ge R_0$ it follows that
\[
\lambda(R)\equiv \left(\int_{0}^{1}\sum_{R-2\le |n|\le R+1}|U(t,n)|^2\right)^{1/2}\ge c e^{-cR \log R}.
\]
\end{theorem}

We remark that this lower bound only uses the fact that the solution is nontrivial and that the constant $c$ in front of the term $R\log R$ only depends on the dimension. 

Theorem \ref{th:DSV} implies Theorem \ref{th:DV}.  The decay conditions at times $t=0$ and $t=1$ imply upper bounds for the term $\lambda(R)$. Indeed, monotonicity results from \cite{JLMP, FBV} show that
\begin{equation}\label{eq:decay}
\|e^{\gamma|n|\log |n|}U(0)\|_{\ell^2(\Z^d)}+\|e^{\gamma|n|\log |n|}U(1)\|_{\ell^2(\Z^d)}<\infty
\end{equation}
for some fixed $\gamma$ implies 
$\|e^{\gamma|n|\log |n|}U(t)\|_{\ell^2(\Z^d)}<\infty$ for all $t\in[0,1]$. Hence, if \eqref{eq:decay} is satisfied,
\[\lambda(R) \le Ce^{-\gamma R\log R}\]
for a positive constant $C$. Thus, by letting $R$ tend to infinity we arrive to a contradiction if $\gamma$ is large enough, since the upper bound decays faster than the lower bound, and therefore $U\equiv0$ if \eqref{eq:decay} is satisfied for $\gamma>\gamma_0$ where $\gamma_0$ depends only on the dimension. However, these results are not sharp.   
We know  that the bound can be improved to
$\exp(-|n|(\log|n|+\mu))$ for some large constant $\mu$. For the free equation ($V=0$), the condition  $\mu>\log 2-1$ implies the uniqueness,  and the question is if for bounded potential the uniqueness result holds with the same range of $\mu$.

Further uniqueness results for solutions of discrete Schr\"odinger type equations, that are inspired by the works of Escauriaza, Kenig, Ponce, and Vega on the continuous case,  can be found in \cite{AR, LM, FBJ, FBGI}.


\end{document}